\documentclass{article}
\usepackage{amsmath,amssymb,amsthm}
\usepackage{graphics}

\newtheorem{theorem}{Theorem}
\newtheorem{remark}{Remark}

\title {Spherical two-distance sets}

\author {Oleg R. Musin \thanks{Department of Mathematics, University of Texas at Brownsville.   oleg.musin@utb.edu}}

\begin{document}
\date{}
\maketitle

\begin{abstract} A set $S$ of unit vectors in $n$-dimensional Euclidean space is called spherical two-distance set, if there are  two numbers $a$ and $b$ so that the inner products  of  distinct vectors of $S$ are either $a$ or $b$. It is known that the largest  cardinality $g(n)$ of  spherical two-distance sets does not exceed $n(n+3)/2$. This upper bound is known to be tight for $n=2,6,22$. The set of mid-points of the edges of a regular simplex gives   the lower bound $L(n)=n(n+1)/2$ for $g(n)$.

In this paper using the so-called  polynomial method it is proved that  for nonnegative $a+b$ the largest cardinality of $S$ is not greater than $L(n)$. For the case $a+b<0$ we propose upper bounds on $|S|$ which are based on Delsarte's method.
 Using this we show that $g(n)=L(n)$ for $6<n<22,\; 23<n<40$, and $g(23)=276$ or 277.
\end {abstract}

\section{Introduction}

A set $S$ in Euclidean space ${\bf R}^n$ is called a {\it two-distance set}, if there are two distances $c$ and $d$, and the distances between pairs of points of $S$ are either $c$ or $d$.
If a two-distance set $S$ lies in the unit sphere ${\bf S}^{n-1}$, then $S$ is called {\it spherical two-distance set.} In other words, $S$ is a set of unit vectors, there are  two real numbers $a$ and $b$,  $-1\le a,b<1$, and inner products  of  distinct vectors of $S$ are either $a$ or $b$.

The ratios of distances of two-distance sets are quite restrictive. Namely, Larman, Rogers, and Seidel \cite {LRS} have proved the following fact:{ if the cardinality of a two-distance set $S$ in ${\bf R}^n,\;$  with distances $c$ and $d, \; c< d,$ is greater than  $2n+3$,  then the ratio $c^2/d^{\,2}$ equals $(k-1)/k$ for an integer   $k$ with
$$
2\le k \le \frac{1+\sqrt{2n}}{2}\,.
$$
}

Einhorn and Schoenberg \cite{ES} proved that there are finitely many two-distance sets $S$ in ${\bf R}^n$ with cardinality $|S|\ge n+2$. Delsarte, Goethals, and Seidel \cite{Del2} proved that the largest cardinality of spherical two-distance sets in ${\bf R}^n$ (we denote it by $g(n)$) is bounded by $n(n+3)/2$, i.e.,
$$
g(n)\le \frac{n(n+3)}{2}\,.
$$
 Moreover, they give examples of spherical two-distance sets  with $n(n+3)/2$ points for $n=2, 6, 22$. {(Therefore, in these dimensions we have equality $g(n)=n(n+3)/2$.)  }
Blockhuis \cite{Blo1} showed that the cardinality of (Euclidean) two-distance sets in ${\bf R}^n$ does not exceed $(n+1)(n+2)/2$.

 The standard unit vectors $e_1,\ldots,e_{n+1}$ form an orthogonal basis of ${\bf R}^{n+1}$. Denote by $\Delta_n$ the  regular simplex with vertices $2e_1,\ldots,2e_{n+1}$. Let $\Lambda_n$ be the set of points $e_i+e_j, \; 1\le i<j\le n+1.$ Since $\Lambda_n$ lies in the hyperplane $\sum {x_k}=2$, we see that $\Lambda_n$ represents a spherical two-distance set in ${\bf R}^n$. On the other hand, $\Lambda_n$ is the set of mid-points of the edges of $\Delta_n$. Thus,
$$
g(n)\ge |\Lambda_n|=\frac{n(n+1)}{2}\,.
$$

For $n<7$ the largest cardinality of Euclidean two-distance sets is $g(n)$, where $g(2)=5, \; g(3)=6, \; g(4)=10, \; g(5)=16$, and $g(6)=27$ (see \cite{Lis}). However, for $n=7,8\; $ Lison\v{e}k \cite{Lis} discovered non-spherical maximal two-distance sets of the cardinality 29 and 45 respectively.

In this paper we prove that
$$
g(n)=\frac{n(n+1)}{2}\,, \; \mbox{ where } \; 6<n<40, \; n\ne 22, 23,
$$
and $g(23)=276$ or 277.
This proof (Section 4) is based on the new sharp upper bound ${n+1\choose 2}$  for spherical two-distance sets with $a+b\ge0$ (Section 2), and on the Delsarte bounds for spherical two-distance sets in the case $a+b<0$. 

\section{Linearly independent polynomials }

The upper bound $n(n+3)/2$ for spherical two-distance sets \cite{Del2}, the bound ${n+2\choose 2}$
for Euclidean two-distance sets \cite{Blo1}, as well as the bound ${n+s\choose s}$  for $s-$distance sets
 \cite{BBS,Blo2}  were obtained by the  polynomial method. The main idea of this method is the following: to associate sets to  polynomials and show that these polynomials are linearly independent as  members of the corresponding vector space.

 Now we apply this idea to improve upper bounds for spherical two-distance sets with $a+b\ge 0.$

\begin{theorem} Let $S$ be a spherical two-distance set in ${\bf R}^n$ with inner products $a$ and $b$. If $a+b\ge0$, then
$$
|S|\le \frac{n(n+1)}{2}\,.
$$
\end{theorem}
\begin{proof} Let
$$
F(t):=\frac{(t-a)(t-b)}{(1-a)(1-b)}\, .
$$
For a unit vector $y\in {\bf R}^n$ we define the function $F_y:{\bf S}^{n-1}\to {\bf R}$ by
$$
F_y(x):=F(\langle x, y\rangle),\quad x\in {\bf R}^{n}, \; ||x||=1.
$$

Let $S=\{x_1,\ldots,x_m\}$ be an $m$-element set. Denote $f_i(x):=F_{x_i}(x)$. Since
$$
f_i(x_j)=\delta_{i,j} \, , \eqno (1)
$$
the quadratic polynomials $f_i,\; i=1,\ldots,m,$ are linearly independent.

Let $e_1, \ldots, e_n$ be a basis of ${\bf R}^n$. Let $L_i(x):=\langle x, e_i\rangle, \; x\in {\bf S}^{n-1}$. Then the  linear polynomials $L_1,\ldots, L_n$ are also linearly independent.

Now we show that if $a+b\ge0$, then $f_1,\ldots, f_m, L_1,\ldots, L_n$ form a linearly independent system of polynomials. 
Assume the converse. Then
$$
\sum\limits_{k=1}^n {d_kL_k(x)}=\sum\limits_{i=1}^m {c_if_i(x)},
$$
where there are nonzero $d_k$ and $c_i$.

Let
$$v=d_1e_1+\ldots+d_ne_n.$$
Then
$$
\langle x, v\rangle=\sum\limits_{i} {c_if_i(x)}. \eqno (2)
$$
For $x=x_i$ in $(2)$, using $(1)$, we get
$$
c_i=\langle x_i, v\rangle.
$$
Take $x=v$ and $x=-v$ in $(2)$. Then we have
$$
||v||^2=\sum\limits_{i} {c_if_i(v)}=\sum\limits_{i} {c_iF(c_i)}, \eqno (3)
$$
$$
-||v||^2=\sum\limits_{i} {c_if_i(-v)}=\sum\limits_{i} {c_iF(-c_i)}. \eqno (4)
$$
Subtracting $(3)$ from $(4)$, we obtain
$$
-||v||^2=\frac{a+b}{(1-a)(1-b)}\,\sum\limits_i{c_i^2}.
$$
This yields $v=0$,  a contradiction.

Note that the dimension of the vector space of quadratic polynomials on the sphere ${\bf S}^{n-1}$ is $n(n+3)/2$. Therefore,
$$
\dim{\{f_1,\ldots,f_m,L_1,\ldots,L_n\}}=m+n\le \frac{n(n+3)}{2}.
$$
Thus, $|S|=m\le n(n+1)/2$.
\end{proof}

Denote by $\rho(n)$ the largest possible cardinality of spherical two-distance sets in ${\bf R}^{n}$ with $a+b\ge0.$

\begin{theorem} If $n\ge7$, then
$$
\rho(n)=\frac{n(n+1)}{2}.
$$
\end{theorem}
\begin{proof} Theorem 1 implies that $\rho(n)\le n(n+1)/2$. On the other hand, the set of mid-points of the edges of a regular simplex has $n(n+1)/2$ points and $a+b\ge0$ for $n\ge7$. Indeed, for this spherical two-distance set we have
$$
a=\frac{n-3}{2(n-1)}\,, \quad b=\frac{-2}{n-1}\,.
$$
Thus,
$$
a+b=\frac{n-7}{2(n-1)}\ge0.
$$
\end{proof}

\section{Delsarte's method  for two-distance sets}

Delsarte's method is widely used in coding theory and discrete geometry for finding bounds for error-correcting codes, spherical codes, and sphere packings (see \cite{CS, Del2, Kab}). In this method upper bounds for spherical codes are given by the following theorem:

\begin{theorem}[\cite{Del2,Kab}]
Let $T$ be a subset of the interval $[-1,1]$. Let $S$ be a set of unit vectors in ${\bf R}^n$ such that the set of inner products  of  distinct vectors of $S$ lies in $T$. Suppose a polynomial $f$ is a nonnegative linear combination of Gegenbauer polynomials $G_k^{(n)}(t)$, i.e.,
$$f(t)=\sum\limits_k {f_kG_k^{(n)}(t)}, \; \mbox{ where }\; f_k\ge 0.$$
If $f(t)\le0$ for all $t\in T$ and $f_0>0$, then
$$ |S|\le \left\lfloor\frac{f(1)}{f_0}\right\rfloor$$
\end{theorem}

There are many ways to define Gegenbauer (or ultraspherical) polynomials
$G_k^{(n)}(t)$.  $G_k^{(n)}$ are a special case of  Jacobi polynomials
$P_k^{(\alpha,\beta)}$ with $\alpha=\beta=(n-3)/2$ and with normalization  $G_k^{(n)}(1)=1$.
Another way to define
$G_k^{(n)}$ is the recurrence formula:
$$ G_0^{(n)}=1,\; G_1^{(n)}=t,\; \ldots,\; G_k^{(n)}=\frac {(2k+n-4)\,t\,G_{k-1}^{(n)}-(k-1)\,G_{k-2}^{(n)}} {k+n-3}.$$
For instance,
$$
G_2^{(n)}(t)=\frac{nt^2-1}{n-1},
$$
$$
G_3^{(n)}(t)=\frac{(n+2)t^3-3t}{n-1},
$$
$$
G_4^{(n)}(t)=\frac{(n+2)(n+4)t^4-6(n+2)t^2+3}{n^2-1}.
$$

Now for given  $n, a, b$ we introduce  polynomials $P_i(t),\; i=1,\ldots,5$.

\medskip

\noindent{$i=1: \;
P_1(t)=(t-a)(t-b)=f_0^{(1)}+f_1^{(1)}t+f_2^{(1)}G_2^{(n)}(t).$

\medskip

\noindent{$i=2: \;
P_2(t)=(t-a)(t-b)(t+c)=f_0^{(2)}+f_1^{(2)}t+f_2^{(2)}G_2^{(n)}(t)+f_3^{(2)}G_3^{(n)}(t)$, where $c$ is defined by the equation $f_1^{(2)}=0$.

\medskip

\noindent{$i=3: \;
P_3(t)=(t-a)(t-b)(t+a+b)=f_0^{(3)}+f_1^{(3)}t+f_2^{(3)}G_2^{(n)}(t)+f_3^{(3)}G_3^{(n)}(t).$ Note that $f_2^{(3)}=0.$

\medskip

\noindent{$i=4: \;
P_4(t)=(t-a)(t-b)(t^2+c\,t+d)=\sum {f_k^{(4)}G_k^{(n)}(t)}$, where $c$ and $d$ are defined by the equations $f_1^{(4)}=f_2^{(4)}=0$.

\medskip

\noindent{$i=5: \;
P_5(t)=(t-a)(t-b)(t^2+c\,t+d)=\sum {f_k^{(5)}G_k^{(n)}(t)}$, where $c$ and $d$ are defined by the equations $f_2^{(5)}=f_3^{(5)}=0$.

\medskip

Denote by $D_i^{(n)}$ the set of all pairs $(a,b)$ such that the  polynomial $P_i(t)$ is well defined, all $f_k^{(i)}\ge 0$, and $f_0^{(i)}>0.$
For instance,
$$
D_1^{(n)}=\left\{(a,b)\in {\bf I}^2: f_1^{(1)}=-a-b\ge0,\; f_0^{(1)} = ab+\frac{1}{n}>0\right\},
$$
$$
D_2^{(n)}=\left\{(a,b)\in {\bf I}^2: a+b\ne0,\; c\ge a+b,\;
f_0^{(2)} = abc+\frac{c-a-b}{n}>0\right\},
$$
where
$$
{\bf I}=[-1,1), \quad c=\frac{ab(n+2)+3}{(n+2)(a+b)}.
$$

Let
$$
U_i^{(n)}(a,b):=\frac{P_i(1)}{f_0^{(i)}}.
$$

Note that we have $P_i(a) = P_i(b) = 0$. Then Theorem 3 yields

\begin{theorem} Let $S$ be a spherical two-distance set in ${\bf R}^n$ with inner products $a$ and $b$. Suppose $(a,b)\in D_i^{(n)}$ for some $i,\; 1\le i\le 5$. Then
$$
|S|\le U_i^{(n)}(a,b).
$$
\end{theorem}

Let $S$ be a spherical two-distance set in ${\bf R}^n$ with inner products $a$ and $b$, where $a> b$. Let $c=\sqrt{2-2a}, \; d=\sqrt{2-2b}$. Then $c$ and $d$ are the Euclidean distances of $S$.

Let
$$
b_k(a)=\frac{ka-1}{k-1}\,.
$$

If $k$ is defined by the equation: $b_k(a)=b$, then $(k-1)/k=c^2/d^{\,2}$.
Therefore,
if $|S|>2n+3$,   then $k$ is an integer number and $k\in\{2,\ldots,K(n)\}$ \cite{LRS}. Here, $K(n) = \lfloor \frac{1+\sqrt{2n}}{2} \rfloor$.


Denote by $D_{i,k}^{(n)}$ the set of all real numbers $a$ such that
$
(a,b_k(a))\in D_i^{(n)}.
$
Let
$$ R_{i,k}^{(n)}(a):=\left\{
\begin{array}{l}
U_i^{(n)}(a,b_k(a)) \; \mbox { for } \; a\in D_{i,k}^{(n)}\\
\infty  \qquad  \qquad \quad  \; \mbox { for } \;  a\notin D_{i,k}^{(n)}
\end{array}
\right.
$$
$$
Q_k^{(n)}(a):=\min\limits_i \left\{R_{i,k}^{(n)}(a)\right\}
$$
Then Theorem 4 yields the following bound for $|S|$:
\begin{theorem} Let $S$ be a spherical two-distance set in ${\bf R}^n$ with inner products $a$ and
$b_k(a)$.  Then
$$
|S|\le Q_k^{(n)}(a).
$$
\end{theorem}

\medskip

\begin{figure}
\includegraphics{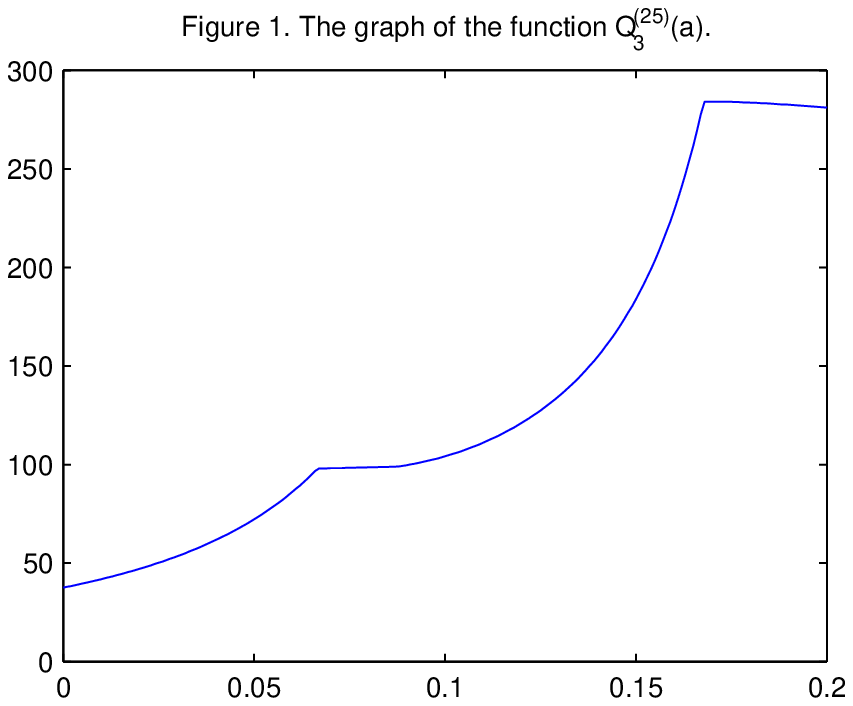}
\end{figure}

\medskip

\medskip

Consider the case  $a+b_k(a)<0$. Since $b_k(a)\ge-1$, we have
$$
a\in I_k:=\left[\frac{2-k}{k}\,, \frac{1}{2k-1}\right).
$$

\begin{remark}
Actually, the polynomials $P_i$ are chosen such that the maximum of  $Q_k^{(n)}(a)$ on $I_k$ minimize the Delsarte bound (Theorem 3).
Clearly, $Q_k^{(n)}(a)$ is a piecewise rational function on $I_k$. It is not hard to find explicit expressions for $Q_k^{(n)}(a)$ and  to compute its maximum on $I_k$ numerically.  For instance,
$\max{\{Q_3^{(25)}(a): a\in I_3=[-1/3,1/5)\}} \approx 284.14$ (see Fig. 1).


\end{remark}

\section{Maximal spherical two-distance sets}

In this section we use Theorem 5 to bound the cardinality of a spherical two-distance set with $a + b < 0$.

  Let $S$,  $|S|>2n+3$, be a spherical two-distance set in ${\bf R}^n$ with inner products $a$ and
$b_k(a)$. Then $k\in\{2,\ldots,K(n)\}$, and $-1\le b_k(a)<a<1$.

Let $\tilde K(n):=\max\{K(n),2\}$. For  given $n$ and $k=2,\ldots,\tilde K(n)$,  we denote by $\Omega(n,k)$ the set of all spherical two-distance sets $S$ in ${\bf R}^n$  with 
$a+b_k(a)<0$. Denote by $\omega(n,k)$ the largest cardinality of $S\in \Omega(n,k)$.

Let
$$
 \varphi(n,k):=\sup\limits_{a\in I_k}\left\{Q_k^{(n)}(a)\right\},
$$
$$
\widehat\omega(n,k):= \max\{\lfloor\varphi(n,k)\rfloor,\, 2n+3\}.
$$

Let us denote by $\widehat\omega(n)$ the maximum of numbers $\widehat\omega(n,2),\ldots,\widehat\omega(n,\tilde K(n))$, and by $\omega(n)$ we denote the largest cardinality of a two-distance set $S$ in ${\bf S}^{n-1}$ with  
$a+b<0$. Then $g(n)=\max\{\omega(n),\rho(n)\}$.

 Since Theorem 5 implies
$\omega(n,k)\le\widehat\omega(n,k),$
we have

\begin{theorem}
 $g(n)\le \max\{\widehat\omega(n),\rho(n)\}$.
\end{theorem}

Finally,  for $g(n)$  we have the following bounds: $\rho(n)\le g(n)\le \max\{\widehat\omega(n),\rho(n)\}$.
Recall that $\rho(n)=n(n+1)/2$ for $n\ge7$. For $\widehat\omega(n), \; 7\le n\le 40,$
we obtain the computational results gathered in Table 1.

\medskip

Since  $\widehat\omega(n)\le \rho(n)$ for  $6<n<40, \; n\ne 22,23$, for these cases we  have  $g(n)=\rho(n)$. For $n=23$ we obtain $g(23)\le 277$. But $g(23)\ge\rho(23)=276$. This proves the following theorem:

\begin{theorem} If $\; 6<n<22\; $ or $\; 23<n<40$, then
$$
g(n)=\frac{n(n+1)}{2}\,.
$$
For $n=23$ we have
$$
g(23)=276 \; \mbox{ or } \; 277.
$$
\end{theorem}

\medskip

\begin{remark}
The case $n=23$ is very interesting. In this dimension the maximal number of equiangular lines (or equivalently, the maximal cardinality of a two-distance set with $a+b=0$) is $276$ \cite{LeS}. On the other hand, $|\Lambda_{23}|=276$. So in $23$ dimensions we have  two very different two-distance sets with $276$ points.

Note that $\max{\{Q_3^{(23)}(a):\; a\in I_3\}}\approx 277.095$. So this numerical bound is not far from $277$. Perhaps stronger tools, such as semidefinite programming bounds, are needed here to prove that $g(23)=276.$
\end{remark}

\begin{remark}
Our numerical calculations show that the barrier $n=40$ is in fact fundamental: LP bounds are incapable of resolving the $n\ge40,\; k=2$ case. That means a new idea is required to deal with $n\ge40.$
\end{remark}

\medskip

\medskip

{\centerline{
{ Table 1}. $\widehat\omega(n)$ and $\rho(n)$. The last column gives the $k$ with  $ \widehat\omega(n)=\widehat\omega(n,k).$
}}

\medskip

\centerline {
\begin{tabular}{|r|r|r|r|}
\hline
$n$ & $\widehat\omega$ & $\rho$ & $k$ \\
\hline
7&28&28&2\\
8  &  31 &   36 &    2\\
     9  &  34 &   45 &    2\\
    10  &   37  &   55  &    2\\
    11  &   40  &   66   &   2\\
    12  &   44  &   78   &   2\\
    13  &   47  &   91   &   2\\
    14  &  52 &  105  &   2\\
    15  &  56 &  120  &   2\\
    16  &  61  & 136  &   2\\
    17   & 66  & 153 &    2\\
    18  &  76 &  171  &   3\\
    19  &  96 &  190   &  3\\
    20  & 126  & 210  &   3\\
    21  & 176 &  231  &   3\\
    22  & 275  & 253  &   3\\
    23 &  277 &  276  &   3\\
    24  & 280  & 300 &    3\\
    25  & 284 &  325   &  3\\
    26  & 288  & 351  &   3\\
    27  & 294  & 378   &  3\\
    28  & 299  & 406  &   3\\
    29 &  305  & 435  &   3\\
    30  & 312 &  465   &  3\\
    31  & 319 &  496   &  3\\
    32  & 327 &  528  &   3\\
    33  & 334 &  561  &   3\\
    34 &  342 &  595   &  3\\
    35  & 360  & 630   &  2\\
    36  & 416  & 666  &   2\\
    37 &  488  & 703  &   2\\
    38  & 584  & 741  &   2\\
    39  & 721 &  780  &   2\\
    40 &  928 &  820  &   2 \\
\hline
\end{tabular}
}

\medskip

\medskip

\begin{remark} It is known that for  $n=3,7,23$ maximal spherical two-distance sets are not unique, and for $n=2,6, 22$, when $g(n)=n(n+3)/2$, these sets are unique up to isometry. Lison\v{e}k \cite{Lis} confirmed the maximality and uniqueness of previously known sets for $n=4,5,6$. For all other $n$ the problem of uniqueness of maximal two-distance sets is open. We think that for $7<n<46, \; n\ne22,23$  maximal spherical two-distance sets in ${\bf R}^n$ are unique and congruent to $\Lambda_n$.

\end{remark}

\medskip

\medskip

\medskip

\medskip

\noindent {\bf\Large Acknowledgments}

\medskip

\medskip


The author thanks Alexander Barg and Frank Vallentin for useful discussions and comments.


\medskip

\medskip




\begin{thebibliography}{99}

\bibitem{BBS} E. Bannai, E. Bannai, and D. Stanton, An upper bound for the cardinality of an $s$-distance set in real Euclidean space, {\it Combinatorica}, {\bf 3} (1983), 147-152.

\bibitem{Blo1} A. Blokhuis,  A new upper bound for the cardinality of  $2$-distance set in  Euclidean space, {\it Ann. Discrete Math.},  {\bf 20}  (1984), 65-66.

\bibitem{Blo2} A. Blokhuis, Few-distance sets, CWI Tract 7 (1984).

\bibitem{CS}
J. H. Conway and N. J. A. Sloane, Sphere Packings, Lattices, and Groups,  Springer-Verlag, New York-Berlin, 1988.

\bibitem{Del2}
Ph. Delsarte, J. M. Goethals and J. J. Seidel, Spherical codes and designs, {\it Geom. Dedic.}, {\bf 6}, 1977, 363-388.

\bibitem{ES}
S. J. Einhorn and I. J. Schoenberg, On Euclidean sets having only two distances between points I, II, {\it Indag. Math.}, {\bf 28} (1966), 479-488, 489-504. (Nederl. Acad. Wetensch. Proc. Ser. A69)

\bibitem{Kab}
G. A. Kabatiansky and V. I. Levenshtein, Bounds for packings on a sphere and in space,
{\it Problems of Information Transmission}, {\bf 14}(1), 1978, 1-17.

\bibitem{LRS}
D. G. Larman, C. A. Rogers, and J. J. Seidel, On two-distance sets in Euclidean space, {\it Bull. London Math. Soc.}, {\bf 9} (1977), 261-267.

\bibitem{LeS}
P. W. H. Lemmens and J. J. Seidel, Equiangular lines, {\it J. Algebra}, {\bf 24} (1973), 494-512.

\bibitem{Lis}
P. Lison\v{e}k, New maximal two-distance sets, {\it J. Comb. Theory, Ser. A}, {\bf 77} (1997), 318-338.

 \end{thebibliography}
\end{document}